\documentclass[11pt,fleqn,twoside]{article}
\usepackage{amsfonts,amssymb,latexsym}
\makeatletter
\newcommand{\prava}[1]{\small\it
\begin{flushleft}
Copyright \copyright \ 1999 by  #1
\end{flushleft}}

\newcommand{\name}[1]{\begin{flushleft}
                       \LARGE \bf #1
                       \end{flushleft}\vspace{-3mm}}

\newcommand{\Author}[1]{\begin{flushleft}
                       \it #1 \end{flushleft}}

\newcommand{\Adress}[1]{\begin{flushleft}
                       \it #1 \end{flushleft}}

\newcommand{\Date}[1]{\begin{flushleft}
                      \small  \it #1 \end{flushleft}}

\newcommand{\ehkol}{Author \ name}
\newcommand{\ohkol}{Article \ name}
\renewcommand{\@evenhead}{
\hspace*{-3pt}\raisebox{-15pt}[\headheight][0pt]{\vbox{\hbox to \textwidth 
{\thepage \hfil \ehkol}\vskip4pt \hrule}}}
\renewcommand{\@oddhead}{
\hspace*{-3pt}\raisebox{-15pt}[\headheight][0pt]{\vbox{\hbox to \textwidth 
{\ohkol \hfil \thepage}\vskip4pt\hrule}}}
\renewcommand{\@evenfoot}{}
\renewcommand{\@oddfoot}{}

\newcommand{\be}{\begin{equation}}
\newcommand{\ee}{\end{equation}}
\newcommand{\ba}{\hspace*{-5pt}\begin{array}}
\newcommand{\ea}{\end{array}}

\newcommand{\ds}{\displaystyle}
\makeatother

\newcommand{\symb}{\vartriangleright\!\! <}

\begin{document}

\thispagestyle{empty}
\setcounter{page}{222}

\renewcommand{\ehkol}{B.A. Kupershmidt}
\renewcommand{\ohkol}{On the Nature of the Virasoro Algebra}

\begin{flushleft}
\footnotesize \sf
Journal of Nonlinear Mathematical Physics \qquad 1999, V.6, N~2,
\pageref{kupershmidt_4-fp}--\pageref{kupershmidt_4-lp}.
\hfill {\sc Article}
\end{flushleft}

\vspace{-5mm}

\renewcommand{\footnoterule}{}
{\renewcommand{\thefootnote}{} 
 \footnote{\prava{B.A. Kupershmidt}}}

\name{On the Nature of the Virasoro Algebra}\label{kupershmidt_4-fp}

\Author{Boris A. KUPERSHMIDT}

\Adress{The University of Tennessee Space Institute,
Tullahoma, TN 37388, USA\\[1mm]
E-mail: bkupersh@utsi.edu}

\Date{Received November 24, 1998; Accepted December 24, 1998}

\begin{flushright}
\begin{minipage}{7cm}
\small  \bfseries \itshape
To the great algebraist Victor Kac, on occasion of his
55\,$^{th}$ birthday.
\end{minipage}
\end{flushright}

\begin{abstract}
\noindent
The multiplication in the Virasoro algebra
\[
[e_p, e_q] = (p - q) e_{p+q} + \theta \left(p^3 - p\right)
\delta_{p + q}, \qquad p, q \in {\mathbf Z},
\]
\[
[\theta, e_p] = 0,
\]
comes from the commutator $[e_p, e_q] = e_p * e_q - e_q * e_p$ in a
quasiassociative algebra with the multiplication
\renewcommand{\theequation}{$*$}
\be
\ba{l}
\ds e_p * e_q = - {q (1 + \epsilon q) \over 1 + \epsilon (p + q)} e_{p+q}
+ {1 \over 2} \theta \left[p^3 - p + \left(\epsilon - \epsilon^{-1}
\right) p^2 \right]
\delta^0_{p+q},
\vspace{3mm}\\
\ds e_p * \theta = \theta* e_p = 0.
\ea
\ee
The multiplication in a quasiassociative algebra ${\cal R}$ satisf\/ies
the property
\renewcommand{\theequation}{$**$}
\be
a * (b * c) - (a * b) * c = b * (a * c) - (b * a) * c, \qquad
a, b, c \in {\cal R}.
\ee
This property is necessary and suf\/f\/icient for the Lie algebra
{\it Lie}$({\cal R})$ to have a phase space.
The above formulae are put into a cohomological framework,
with the relevant complex being dif\/ferent from the Hochschild
one even when the relevant quasiassociative algebra ${\cal R}$ becomes
associative.  Formula $(*)$ above also has a dif\/ferential-variational
counterpart.
\end{abstract}

\section{Introduction}

Quasiassociative algebras, originally discovered by Vinberg [8]--[10]
and Koszul~[3] in the 1960's in the study of homogeneous convex cones,
appear also as an underlying structure of those Lie algerbras that possess
a phase space. Namely, for a given Lie algebra ${\cal G}$, the following
three conditions are equivalent [5]:

\begin{enumerate}
\item[(i)]  ${\cal G} = Lie ({\cal R})$ for some quasiassociative
algebra ${\cal R}$;

\item[(ii)]  Let $\rho: {\cal G} \rightarrow \mbox{End}({\cal{G}}^*)$
be a representation, not necessarily coadjoint one,
such that on the semidirect sum Lie algebra
${\cal G} {\mathop {\symb}\limits_{\rho}} {\cal G}^*$,
the symplectic form is a 2-cocycle;

\item[(iii)] The natural Poisson bracket on the Lie algebra
${\cal G} {\mathop {\symb}\limits_{\rho}} {\cal G}^*$
is compatible with the canoni\-cal Poisson bracket.

\end{enumerate}

Thus, the quasiassociative algebras form a natural category from
the point of view of Classical and Quantum mechanics. A list of Lie algebras
with a phase space, given in~[5], includes such non-evident cases as Lie
algebras of vector f\/ields on ${\mathbf R}^n$ and current algebras.
One of the principal Lie algebras of physical interest, the Virasoro algebra,
has been, however, conspiciously under-privileged so far.
Its underlying quasiassociative structure is treated in the next two Sections.
(Still more Lie algebras with a phase space can be found in Chapter~2 in~[7].)

Before leaving the phase-space perspective for more mathematical matters,
let me make two comments. First, the category of quasiassociative algebras
is closed with respect to the operation of phase-space extension,
unlike the smaller category of associative algebras: if ${\cal R}$
is quasiassociative then so is T$^*{\cal R}$, where multiplication in
T$^*{\cal R}$ is given by the formula~[5]
\renewcommand{\theequation}{\arabic{section}.\arabic{equation}{\rm a}}
\setcounter{equation}{0}
\be
\left(\begin{array}{c} x \\ \bar x
\end{array} \right) * \left(\begin{array}{c} y \\ \bar y \end{array}
\right) = \left(\begin{array}{c} x * y \\ x * \bar y \end{array} \right),
\qquad x, y \in {\cal R}, \quad \bar x, \bar y \in {\cal R}^* = \mbox{Hom}
({\cal R}, \ldots ),
\ee
\renewcommand{\theequation}{\arabic{section}.\arabic{equation}{\rm b}}
\setcounter{equation}{0}
\be
\langle x * \bar y, y \rangle  = - \langle  \bar y, x * y \rangle .
\ee
Second, if $\rho: {\cal{G}} \rightarrow \mbox{End} ({\cal G}^*)$
is the representation staring in the properties (ii) and (iii) above,
then the associated quasiassociative multiplication on ${\cal G}$
is given by the formula
\renewcommand{\theequation}{\arabic{section}.\arabic{equation}}
\setcounter{equation}{1}
\be
x * y = \rho^d (x) (y),
\ee
where $\rho^d: {\cal G} \rightarrow \mbox{Eng} ({\cal G})$
is the representation dual to $\rho$. The condition for the symplectic form
on ${\cal G} {\mathop {\symb}\limits_{\rho}} {\cal G}^*$
to be a 2-cocycle is then equivalent to the property
\be
\rho^d (x) (y) - \rho^d (y) (x) = [x, y], \qquad
\forall \; x, y \in {\cal G}.
\ee
Thus,
\be
Lie ({\mathrm T}^* {\cal R}) = {\mathrm T}^* Lie ({\cal R}).
\ee
The equation (1.3) appears also in a very dif\/ferent context, as the
condition for the complex of dif\/ferential forms on the Universal Enveloping
Algebra $U({\cal G})$ to be ghost-free (see [6], equations (7.4) and
(7.5).)

Turning back to the Virasoro algebra, we see from formula
$(*)$ that we have what appears to be a central extension of the
corresponding centerless quasiassociative multiplication
\be
e_p * e_q = - {q (1 + \epsilon q) \over 1 + \epsilon (p + q)} e_{p + q},
\qquad  p, q \in {\mathbf Z},
\ee
where $\epsilon$ can be treated as either a formal parameter or a number
such that $\epsilon^{-1} \bar \in {\mathbf Z}$. The next Section contains
a quick verif\/ication that formula (1.5) satisf\/ies the quasiassociativity
property $(**)$. Section 3 is devoted to central extensions of
quasiassociative algebras in general and the algebra (1.5) in particular,
resulting in the formula $(*)$ from the Abstract. In Section~4
we re-interpret in the language of 2-cocycles the property of a bilinear
form to provide a central extension of a quasiassociative algebra;
this interpretation then leads to a complex on the space of cochains
$C^n = \mbox{Hom} ({\cal R}^{\otimes n},\cdot)$. Section~5
generalizes this complex to the case $C^n = \mbox{Hom} ({\cal R}^{\otimes n},
{\cal M})$, where ${\cal R}$ acts nontrivally on ${\cal M}$.
Section~6 deals with the dual objects, homology. The last Section~7 is
devoted to dif\/ferential-variational versions of the
preceding results, for the case when the centerless
Virasoro algebra is replaced by the Lie algebra of vector f\/ields on the
circle with the commutator
\be
[X, Y] = X Y' - X'Y, \qquad {}' = d/d z,
\ee
and the central extension is given by the Gelfand-Fuks 2-cocycle
\be
\omega (X, Y) = \int X Y'''\; d z.
\ee

Appendix~1 contains a short proof that the Virasoro algebra does not
come from an associative one.  Semi-direct sums of quasiassociative
algebras are treated in Appendix~2. In Appendix~3 we prove that if $G$
is a connected Lie group whose Lie algebra ${\cal G}$ comes out of a
quasiassociative algebra then the Lie algebra ${\cal D} (G)$ of vector
f\/ields on $G$ also allows a quasiassociative representation.

\setcounter{equation}{0}

\section{The Centerless Virasoro Algebra}

Suppose a space with a basis $ \{ e_p|p \in G$, a commutative ring$\}$
has the multiplication of the form
\be
e_p * e_q = f(p, q) e_{p+q}.
\ee
Then
\be
\ba{l}
e_p * (e_q * e_r) - (e_p * e_q) * e_r = f(q, r) e_p * e_{q+r}
- f(p, q) e_{p + q} *e_r
\vspace{2mm}\\
\qquad = [f(q,r) f(p,q+r) - f(p,q) f(p+q,r)] e_{p+q+r},
\ea
\ee
so that the quasiassociativity condition $(**)$, the symmetry between
$p$ and $q$, is equivalent to the relation
\be
f (q, r) f(p, q+r) - f(p,q) f(p+q,r) =
 f (p, r) f(q, p+r) - f(q, p) f(p+q,r),
\ee
which can be rewritten as
\be
[f(p,q) - f(q, p)] f(p+q,r) =
f(q, r) f(p, q+r) - f(p, r) f(q, p+r).
\ee

By formula (1.5),
\be
f(p,q) = - {q (1 + \epsilon q) \over 1 + \epsilon (p + q)},
\ee
and we have to check that this $f(p,q)$ satisf\/ies formula (2.4).

First,
\be
f(p, q) - f(q, p) = {1 \over 1 + \epsilon (p+q)} [- q (1 + \epsilon q)
+ p (1 + \epsilon p) ] = p - q,
\ee
so that
\be
e_p * e_q - e_q * e_p = (p - q) e_{p+q},
\ee
guaranteeing that the Lie algebra generated by formula (1.5) is indeed
the centerless Virasoro algebra.

Now, for the LHS of formula (2.4) we obtain
\renewcommand{\theequation}{\arabic{section}.\arabic{equation}$\ell$}
\setcounter{equation}{7}
\be
(p-q) {-r (1 + \epsilon r) \over 1+\epsilon (p + q + r)},
\ee
while for the RHS of formula (2.4) we get
\renewcommand{\theequation}{\arabic{section}.\arabic{equation}$r$}
\setcounter{equation}{7}
\be
\ba{l}
\ds {- r(1 + \epsilon r) \over 1 + \epsilon (q+r)} \cdot {-(q+r) [1 +
\epsilon (q+r)] \over 1 + \epsilon (p + q + r)} - {-r (1 + \epsilon r)
\over 1 + \epsilon (p + r)} \cdot {-(p+r) [1 + \epsilon (p + r)]
\over 1 + \epsilon (p + q + r) }
\vspace{3mm}\\
\ds \qquad =
{- r (1 + \epsilon r) \over 1 + \epsilon (p + q + r) } [ - (q + r) +
(p + r)] = {- r(1 + \epsilon r) \over 1 + \epsilon (p + q+r) }
(p - q),
\ea
\ee
and this is the same as formula (2.8$\ell$).

\medskip

\noindent
{\bf Remark 2.9.} Formula (2.5) is not the only solution of the
equation (2.4) satisfying the Lie boundary condition
\renewcommand{\theequation}{\arabic{section}.\arabic{equation}}
\setcounter{equation}{9}
\be
f (p, q) - f (q, p) = p - q.
\ee
For example,
\be
f (p, q) = \lambda - q, \qquad \lambda = \mbox{const},
\ee
is also a solution. It does not alow a proper central extension,
however.
\setcounter{equation}{0}

\section{Central Extensions of Quasiassociative Algebras}

Let $K$ be a commutative ring over which our quasiassociative algebra
${\cal R}$ is an algebra. Let $\Omega: {\cal R} \times {\cal R}
\rightarrow K$ be a bilinear form. It def\/ines a multiplication on the space
$\tilde {\cal R} = {\cal R} \oplus K$, by the rule
\be
\left(\begin{array}{c} a \\ \alpha \end{array} \right) *
\left(\begin{array}{c} b \\ \beta \end{array} \right) =
 \left(\begin{array}{c} a * b \\ \Omega (a, b) \end{array}\right),
 \qquad a, b \in {\cal R}, \quad  \alpha, \beta \in K.
 \ee

When is $\tilde {\cal R}$ quasiassociative? We have:
\be
\ba{l}
\ds \left(\begin{array}{c}
a \\ \alpha \end{array} \right) * \left( \left(
\begin{array}{c}b \\ \beta \end{array} \right) * \left(
\begin{array}{c} c \\ \gamma \end{array} \right) \right) -
\left( \left(\begin{array}{c} a \\ \alpha \end{array} \right) *
\left(\begin{array}{c} b \\ \beta \end{array} \right) \right) *
\left(\begin{array}{c} c \\ \gamma \end{array} \right)
\vspace{3mm}\\
\ds = \left(\begin{array}{c} a \\ \alpha \end{array} \right) *
\left(\begin{array}{c} b * c \\ \Omega (b,c) \end{array} \right) -
\left(\begin{array}{c} a * b \\ \Omega (a, b) \end{array}\right) *
\left(\begin{array}{c} c \\ \gamma \end{array} \right) =
\left(\begin{array}{c} a * (b * c)  -  (a * b) * c \\
\Omega(a, b * c)  -  \Omega (a * b, c) \end{array}\right).
\ea\hspace{-6.73pt}
\ee
Thus, $\tilde {\cal R}$ is quasiassociative if\/f
\be
\Omega (a, b * c) - \Omega (a * b, c) = \Omega (b, a * c) -
\Omega (b * a, c).
\ee
This can be equivalently rewritten as
\be
\Omega (b, a * c) - \Omega (a, b * c) + \Omega ([a, b], c) = 0,
\ee
where $[a, b] = a * b - b* a$ is the commutator in the Lie algebra
$Lie({\cal R})$. By construction, the bilinear form
\be
\omega (a, b) = \Omega (a, b) - \Omega (b, a)
\ee
def\/ines a central extension of the Lie algebra $Lie({\cal R}$); thus,
$\omega$ is a 2-cocycle on this Lie algebra.

While we are at it, let's look at {\it trivial} 
central extensions of ${\cal R}$. These are produced from the
multiplication 
\be
\left(\begin{array}{c} a \\ \alpha \end{array} \right)
{\mathop{*}\limits_t} \left(\begin{array}{c} b \\
\beta \end{array} \right) =
\left(\begin{array}{c} a * b \\ 0 \end{array} \right)
\ee
by linear transformations of the form
\be
\Phi = \left(\begin{array}{cc} id & 0 \\
\langle u, \cdot \rangle  &  1 \end{array} \right), \qquad
u \in {\cal R}^*.
\ee
Thus, trivial extensions look like
\be
\ba{l}
\left(\begin{array}{c} a \\ \alpha \end{array} \right) *
\left(\begin{array}{c} b \\ \beta \end{array}
\right) = \Phi \left( \Phi^{-1} \left(
\begin{array}{c} a \\ \alpha \end{array} \right)
{\mathop{*}\limits_t}
\Phi^{-1} \left(\begin{array}{c} b \\ \beta \end{array}
\right) \right)
\vspace{3mm}\\
\ds \qquad =\Phi \left( \left(\begin{array}{c} a \\ \cdots \end{array} \right)
{\mathop{*}\limits_t}
\left(\begin{array}{c}b \\ \cdots \end{array} \right) \right) =
 \Phi \left(\begin{array}{c} a * b \\ 0 \end{array} \right) =
 \left(\begin{array}{c} a * b \\ \langle  u, a * b \rangle
 \end{array} \right),
 \ea
 \ee
so that trivial ``2-cocycles'' on ${\cal R}$ are of the form
\be
\Omega (a, b) = \langle u, a * b\rangle , \qquad u \in {\cal R}^*.
\ee
The award of the title ``cocycle'' to $\Omega$ will be justif\/ied
in the next Section, where the criterion (3.4) is recast as
\be
\delta \Omega (a, b, c) = 0.
\ee
Similarly, central extensions dif\/fering by a trivial 2-cocycle are
equivalent: if
\[
\left(\begin{array}{c} a \\ \alpha \end{array} \right)
{\mathop {*}\limits_1}
\left(\begin{array}{c} b \\ \beta \end{array} \right) =
\left(\begin{array}{c} a * b \\ \omega (a, b) + \langle  u, a* b
\rangle \end{array} \right)
\]
and
\[
\left(\begin{array}{c}
a \\ \alpha \end{array} \right)
{\mathop{*}\limits_2}
\left(\begin{array}{c} b \\ \beta \end{array} \right) =
\left(\begin{array}{c} a * b \\ \omega (a, b) \end{array} \right)
\]
are two such extensions, then the transformation $\Phi$ (3.7) takes the
second multiplication into the f\/irst one.


Let us return to the case of the Virasoro algebra. Formula $(*)$ shows
that we have a central extension
\be
\Omega (e_p, e_q) = \varphi (p) \delta^0_{p +q}.
\ee
\be
2 \varphi (p) = p^3 - \epsilon^{-1} p^2 - p + \epsilon p^2 =
- p (1 -\epsilon p) \left(1 + \epsilon^{-1} p\right).
\ee
The condition (3.4), in the notation (2.1) and (3.11), becomes:
\be
\ba{l}
\ds 0 = \Omega (e_q, e_p * e_r) - \Omega (e_p , e_q * e_r) + \Omega
([e_p, e_q], e_r)
\vspace{2mm}\\
\ds \qquad = f (p, r) \Omega (e_q, e_{p + r}) -
f(q, r) \Omega (e_p, e_{q+r}) +(p-q) \Omega (e_{p+q}, e_r)
\vspace{2mm}\\
\ds \qquad = [f (p, r) \varphi (q) - f (q, r) \varphi (p) +
(p - q) \varphi (p +q)] \delta_{p + q + r}^0,
\ea
\ee
which can be rewritten as
\be
(p-q) \varphi (p+q) = \varphi (p) f(q, - p-q) - \varphi (q) f (p, - p-q).
\ee

With $f(p,q)$ and $\varphi (p)$ given by formula (2.5) and (3.12)
respectively, for the $2 \times LHS$ of formula (3.14) we get:
\renewcommand{\theequation}{\arabic{section}.\arabic{equation}$\ell$}
\setcounter{equation}{14}
\be
-(p-q) (p+q) [1 - \epsilon (p+q)] \left[1 + \epsilon^{-1} (p+q)\right],
\ee
while for the $2 \times RHS$ of formula (3.14) we obtain:
\renewcommand{\theequation}{\arabic{section}.\arabic{equation}$r$}
\setcounter{equation}{14}
\be
\ba{l}
\ds - p (1 - \epsilon p) \left(1 + \epsilon^{-1} p\right)
\frac{(p+q) [1-\epsilon (p+q)]}{1 - \epsilon p}
\vspace{3mm}\\
\ds \qquad + q (1 - \epsilon q)
\left(1 + \epsilon^{-1} q\right) \frac{(p+q) [1 - \epsilon (p+q)]}
{1 - \epsilon q}
\vspace{3mm}\\
\ds \qquad = (p+q) [1-\epsilon (p+q)] (q - p) \left[1 +
\epsilon^{-1} (p+q)\right],
\ea
\ee
and this is the same as the expression $(3.15\ell)$.

Thus, we get a central extension of the quasiassociative algebra (2.1),
(2.5). It remains to check that the 2-cocycle $\omega $ (3.5) is indeed
the one entering the Virasoro algebra. We have:
\renewcommand{\theequation}{\arabic{section}.\arabic{equation}}
\setcounter{equation}{15}
\be
\ba{l}
\ds \omega (e_p, e_q) = \Omega (e_p, e_q) - \Omega (e_q, e_p) =
[\varphi (p) - \varphi (q) ] \delta^0_{p+q}
\vspace{2mm}\\
\ds \qquad = {1 \over 2} \left\{ \left[p^3 - p +
\left(\epsilon - \epsilon^{-1} \right) p^2 \right] -
\left[(-p)^3 - (-p) + \left(\epsilon - \epsilon^{-1}\right) (-p)^2\right]
\right\} \delta^0_{p+q}
\vspace{2mm}\\
\ds \qquad = \left(p^3-p\right) \delta^0_{p+q}.
\ea
\ee

\setcounter{equation}{0}

\section{The Quasiassociative Complex}

Let ${\cal M}$ be a $K$-module. Def\/ine the cochains on ${\cal R}$
with values of ${\cal M}$ as
\be
C^n = C^n ({\cal R}, {\cal M}) = \mbox{Hom}_K ({\cal R}^{\otimes n} ,
{\cal M}), \qquad n \in {\mathbf N}; \quad C^0: = {\cal M}.
\ee
In the preceding Section we in ef\/fect met two coboundary operators
$\delta: C^n \rightarrow C^{n+1}$ for $n=1$ and $n=2$, in formulae
(3.9) and (3.4) respectively:
\be
\psi \in C^1 \quad  \Rightarrow \quad
 \delta \psi (a_1, a_2) = \psi (a_1 * a_2),
\ee
\be
\psi \in C^2 \quad \Rightarrow \quad \delta \psi (a_1, a_2, a_3)
= \psi (a_2, a_1 * a_3) - \psi (a_1, a_2 * a_3) + \psi ([a_1, a_2], a_3).
\ee

It is obvious that $\delta^2=0$ on $C^1$, and the roundabout way this
equality was verif\/ied in the preceeding Section actually proves that
\be
H^2({\cal R}) : = H^2 ({\cal R}, K)
\ee
describes the $K$-module of isomorphism classes of 1-dimensional central
extensions of ${\cal R}$ by $K$.

Guided by formulas (4.2) and (4.3), we def\/ine the coboundary operator
$\delta: C^n \rightarrow C^{n+1}$ for all $n \in {\mathbf Z}_+$,
as follows
\be
\delta= 0 \quad  \mbox{on} \quad C^0;
\ee
\[
\psi \in C^n, \quad n \geq 1 \quad \Rightarrow
\]
\renewcommand{\theequation}{\arabic{section}.\arabic{equation}{\rm a}}
\setcounter{equation}{5}
\be
\delta \psi (a_1,\ldots, a_n, a) = \sum^n_{i = 1} (-1)^{i+1} \psi (\ldots
\hat i \ldots, a_i * a) +
\ee
\renewcommand{\theequation}{\arabic{section}.\arabic{equation}{\rm b}}
\setcounter{equation}{5}
\be
\qquad + \sum_{1 \leq i < j \leq n} (-1)^{i+j+1} \psi ([a_i, a_j]\ldots \; \hat i
\ldots \; \hat j \; \ldots \, a).
\ee
The hat over the argument signif\/ies this argument's absence; the last
sum (4.6b) is missing when $n<2$; the right-most argument,
$a=a_{n+1}$, is considered on a dif\/ferent footing from the rest,
$a_1, \ldots, a_n$.  (We see that $H^1 ({\cal R}) = \{\mbox{Ann}
({\cal R} * {\cal R}) \subset {\cal R}^*\}$.)

Before proceeding further, we need to make some minimal skewsymmetry
observations.

\medskip

\noindent
{\bf Def\/inition 4.7.} {\it Suppose $n \geq 3$. If $\kappa$
is such that $2 \leq \kappa < n$, then a cochain $\psi \in C^n$ is
called $\kappa$-skewsymmetric if it is skewsymmetric in its first
$\kappa$ arguments.}

\medskip

\noindent
{\bf Proposition 4.8.} {\it  (i)  For $n=2$, $\delta (C^n)$ is
2-skewsymmetric;\\
(ii)  Suppose $n \geq 3$ and $\psi \in C^n$; if $\psi$ is
$\kappa$-skewsymmetric then so is $\delta \psi$.}

\medskip

\noindent
{\bf Proof.}  (i)  Formula (4.3) makes the claim obvious for
$n=2$;\\
(ii)  For $n \geq 3$, the sums (4.6a) and (4.6b)
each change sign under the transposition $(i, i+1)$ for all $i < \kappa $.
\hfill $\blacksquare$

\medskip

Thereafter we assume that all our cochains are $\kappa$-skewsymmetric for
some f\/ixed $\kappa \geq 2$.

\medskip

\noindent
{\bf Proposition 4.9.} {\it $\delta^2 = 0$ on 2-skewsymmetric cochains.}

\medskip

\noindent
{\bf Proof.} Let $\psi \in C^n, n \geq 2$. Set $\nu = \delta
\psi$:
\[
\hspace*{-8.5pt}
\nu (a_1,\ldots, a_n, z) = \sum^n_{i = 1} (-1)^{i+1} \psi(\ldots \hat i \ldots
, iz) + \sum_{1 \leq i < j \leq n} (-1)^{i+j+1} \psi ([i,j] \ldots \hat i
\ldots \hat j \ldots , z),
\]
where for brevity we write $iz$ insted of $a_i * z$, and $[i, j]$
instead of $[a_i, a_j]$. Then,
\[
\delta \nu (y_1,\ldots, y_{n+1}, t)
\]
\renewcommand{\theequation}{\arabic{section}.\arabic{equation}{\rm a}}
\setcounter{equation}{10}
\be
\qquad \qquad = \sum^{n+1}_{s=1} (-1)^{s+1} \nu (\hat s, st) +
\ee
\renewcommand{\theequation}{\arabic{section}.\arabic{equation}{\rm b}}
\setcounter{equation}{10}
\be
\qquad \qquad
+ \sum_{1 \leq p < q \leq n+1} (-1)^{p+q+1} \nu [(p, q]  \hat p  \hat q,
t),
\ee
where for further bref\/ity we now suppress the ``$\ldots$'' convention.

We shall work out separately the expression (4.11a) and (4.11b).

\noindent
({\it a}) \ We have:
\renewcommand{\theequation}{\arabic{section}.\arabic{equation}}
\setcounter{equation}{11}
\be
\ba{l}
\ds \nu (\hat s,t) = \sum_{i<s} (-1)^{i+1} \psi (\hat i \hat s
, i (st)) + \sum_{i > s} (-1)^i \psi (\hat s \hat i, i (st)) 
\vspace{3mm}\\
\ds \qquad+ \sum_{i < j < s} (-1)^{i+j+1} \psi ([i, j] \hat i \hat j \hat s ,
st) 
\vspace{3mm}\\
\ds \qquad + \sum_{i < s, j>s} (-1)^{i+j} \psi ([i, j] \hat i \hat s
\hat j, st) + \sum_{s< i< j} (-1)^{i+j+1} \psi ([i, j] \hat s \hat i
\hat j, st).
\ea
\ee

Multiplying all this by $(-1)^{s+1}$ and summing on $s$, we get
\[
\{(4.11{\rm a})\} =
\]
\renewcommand{\theequation}{\arabic{section}.\arabic{equation}{\rm a}}
\setcounter{equation}{12}
\be
\qquad  = \sum_{a < b} (-1)^{a+b}\{\psi (\hat a \hat b ,
a (bt)) 
\ee
\renewcommand{\theequation}{\arabic{section}.\arabic{equation}{\rm b}}
\setcounter{equation}{12}
\be
\qquad \qquad  - \psi (\hat a  \hat b, b (at))\}
\ee
\renewcommand{\theequation}{\arabic{section}.\arabic{equation}{\rm c}}
\setcounter{equation}{12}
\be
\qquad + \sum_{a < b< c} (-1)^{a+b+c} \{ \psi ([a, b]
\hat a  \hat b  \hat c, ct) 
\ee
\renewcommand{\theequation}{\arabic{section}.\arabic{equation}{\rm d}}
\setcounter{equation}{12}
\be
\qquad \qquad - \psi ([a,c]  \hat a  \hat b 
\hat c , bt)
\ee
\renewcommand{\theequation}{\arabic{section}.\arabic{equation}{\rm e}}
\setcounter{equation}{12}
\be
\qquad \qquad  + \psi ([b,c]  \hat a  \hat b
\hat c, at)\} ; 
\ee

\noindent
({\it b}) \ We have:
\[
\ba{l}
\ds \nu ([p, q] \hat p \hat q , t) = \psi 
( \hat p \hat q , [p, q] t) +  \sum_{\alpha < p} (-1)^\alpha \psi 
([p, q])  \hat \alpha  \hat p \hat q, \alpha t)
\vspace{3mm}\\
\ds \qquad + \sum_{p < \alpha < q} (-1)^{\alpha + 1} \psi ([p, q]
\hat p \hat \alpha \hat q , \alpha t) + 
\sum_{\alpha > q} (-1)^\alpha \psi ([p, q] \hat p \hat q \hat \alpha,
\alpha t)
\vspace{3mm}\\
\ds \qquad + \sum_{\alpha < p} (-1)^{\alpha + 1} \psi ([[p, q],
\alpha ] \hat \alpha  \hat p  \hat q , t) +
\sum_{p <\alpha < q} (-1)^\alpha \psi ([[ p, q], \alpha]  \hat p 
\hat \alpha  \hat q , t) 
\ea
\]
\[
\ba{l}
\ds \qquad + \sum_{\alpha > q} (-1)^{\alpha + 1} \psi ([[ p,
q], \alpha]  \hat p  \hat q  \hat \alpha, t) 
+ \sum_{i < j < p} (-1)^{i + j+1} \psi ([i, j], [p, q]  \hat i 
 \hat j \hat p  \hat q , t) 
\vspace{3mm}\\
\ds \qquad + \sum_{i< p, p< j< q} (-1)^{i+j} \psi [i, j], [p, q]  
\hat i  \hat p \hat j  \hat q , t) 
+ \sum_{i<p, j>q} (-1)^{i + j+1} \psi ([i, j], [p, q]  \hat i
\hat p  \hat q \hat j , t) 
\vspace{3mm}\\
\ds \qquad + \sum_{p< i<j<q} (-1)^{i+j+1} \psi ([i, j], [p,q] \hat p 
\hat i \hat j \hat q , t) 
+ \sum_{p<i< q<j} (-1)^{i+j} \psi ([i, j], [p, q]  \hat p \hat i 
\hat q  \hat j , t) 
\vspace{3mm}\\
\ds \qquad+ \sum_{q< i<j} (-1)^{i+j+1} \psi ([i, j], [p, q] 
 \hat p  \hat q \hat i  \hat j , t). 
\ea
\]

Multiplying this monstrocity by $(-1)^{p+q+1}$ and summing on
$\{1 \leq p <q \leq n + 1\}$, we f\/ind:
\[
\{(4.11{\rm b})\} =
\]
\renewcommand{\theequation}{\arabic{section}.\arabic{equation}{\rm a}}
\setcounter{equation}{13}
\be
\qquad = - \sum_{a < b} 
(-1)^{a+b} \psi ( \hat a \hat b, [a, b]t) 
\ee
\renewcommand{\theequation}{\arabic{section}.\arabic{equation}{\rm b}}
\setcounter{equation}{13}
\be
\qquad + \sum_{a < b < c} (-1)^{a + b+c} \{ - \psi ([b, c] 
 \hat a \hat b \hat c , at) 
\ee
\renewcommand{\theequation}{\arabic{section}.\arabic{equation}{\rm c}}
\setcounter{equation}{13}
\be
\qquad + \psi ([a, c]  \hat a  \hat b  \hat c , bt) 
\ee
\renewcommand{\theequation}{\arabic{section}.\arabic{equation}{\rm d}}
\setcounter{equation}{13}
\be
\qquad - \psi ([a, b]  \hat a  \hat b  \hat c , ct) 
\ee
\renewcommand{\theequation}{\arabic{section}.\arabic{equation}{\rm e}}
\setcounter{equation}{13}
\be
\qquad + \psi ([[b, c],a] \hat a  \hat b  \hat c , t) 
\ee
\renewcommand{\theequation}{\arabic{section}.\arabic{equation}{\rm f}}
\setcounter{equation}{13}
\be
\qquad - \psi ([[a, c],b] \hat a \hat b \hat c , t) 
\ee
\renewcommand{\theequation}{\arabic{section}.\arabic{equation}{\rm g}}
\setcounter{equation}{13}
\be
\qquad + \psi ([[a, b],c]  \hat a \hat b  \hat c  , t)\} 
\ee
\renewcommand{\theequation}{\arabic{section}.\arabic{equation}{\rm h}}
\setcounter{equation}{13}
\be
\qquad + \sum_{a< b< c< d} (-1)^{a+b+c+d} \{\psi ([a, b], [c, d]
\hat a  \hat b  \hat c  \hat d , t) 
\ee
\renewcommand{\theequation}{\arabic{section}.\arabic{equation}{\rm i}}
\setcounter{equation}{13}
\be
\qquad  - \psi ([a, c], [b, d] \hat a  \hat b  \hat c  \hat d , t) 
\ee
\renewcommand{\theequation}{\arabic{section}.\arabic{equation}{\rm j}}
\setcounter{equation}{13}
\be
\qquad + \psi ([a, d], [b, c] \hat a \hat b \hat c \hat d,t) 
\ee
\renewcommand{\theequation}{\arabic{section}.\arabic{equation}{\rm k}}
\setcounter{equation}{13}
\be
\qquad + \psi ([b, c], [a, d] \hat a \hat b \hat c \hat d ,t) 
\ee
\renewcommand{\theequation}{\arabic{section}.\arabic{equation}{\rm l}}
\setcounter{equation}{13}
\be
\qquad  - \psi ([b, d], [a, c]  \hat a \hat b \hat c \hat d,t) 
\ee
\renewcommand{\theequation}{\arabic{section}.\arabic{equation}{\rm m}}
\setcounter{equation}{13}
\be
\qquad  + \psi ([o, d]  [a, b]  \hat a  \hat b  \hat c  \hat d, t)\}.
\ee

Grouping various terms together, we organize the cancellation scheme as
follows:
\begin{enumerate}
\item[1)]  (4.13a,b) and (4.14a), because of the equality
\renewcommand{\theequation}{\arabic{section}.\arabic{equation}}
\setcounter{equation}{14}
\be
a (bt) - b (at) = [a, b] t 
\ee
being the def\/ining relation $(**)$ of a quasiassociative algebra;

\item[2)] (4.13c) and (4.14d);  (4.13d) and (4.14c);  (4.13e) and
(4.14b); 

\item[3)] (4.14e,f,g) by virtue of the Jacobi identity;

\item[4)] (4.14h,m); (4.14i,l); (4.14j,k); -- all by virtue of $\psi$
being 2-skewsymmetric.  \hfill $\blacksquare$

\end{enumerate}

\renewcommand{\theequation}{\arabic{section}.\arabic{equation}}
\setcounter{equation}{0}

\section{The Quasiassociative Complex with Values in a Module}

In this Section we generalize the coboundary operator
$\delta: C^n\rightarrow C^{n+1}$ given by formula~(4.6), to the case where
the quasiassociative algebra ${\cal R}$ acts nontrivially on ${\cal M}$,
the space where cochains take values.

Suppose $\chi: {\cal R} \rightarrow \mbox{End} ({\cal M})$ is a linear map.
It is natural to call it a representation of ${\cal R}$ if it behaves
the way the left multiplication in ${\cal R}$ does:
\be
\chi (a) \chi (b)- \chi (a * b) = \chi (b) \chi(a) - \chi (b*a), \qquad
\forall \; a, b \in {\cal R}.
\ee
Since this can be rewritten as
\be
[\chi (a), \chi(b)] = \chi ([a, b]),
\ee
we simply have a representation of the underlying Lie algebra $Lie({\cal R})$.
It is interesting that for the purpose of extending the chain complex
(4.6) of the preceeding Section, this natural and proper def\/inition is
insuf\/f\/icient; a stronger one is required. This insuf\/f\/iciency can be
seen as follows.

Let $\psi: {\cal R} \rightarrow {\cal M}$ be a 1-cochain. By formula (4.2),
we should now have
\be
\delta \psi (a, b) = \psi (a * b) + c_1 \chi (a) \psi (b) + c_2 \chi (b)
\psi(a),
\ee
with some constants $c_1$ and $c_2$. If we f\/ix $m \in C^0 = {\cal M}$
and consider the natural def\/inition for the operator
$\delta: C^0 \rightarrow C^1$,
\be
\delta m (a) = \chi (a) (m),
\ee
then
\be
\ba{l}
\left(\delta^2 m\right) (a, b) = \delta m (a * b) + c_1 \chi (a)
\delta m (b) + c_2 \chi (b) \delta m (a)
\vspace{2mm}\\
\qquad = (\chi(a * b) + c_1 \chi (a) \chi (b) + c_2 \chi (b) \chi(a)) (m).
\ea
\ee
This expression has no reasons to vanish unless we change the def\/inition
of (left) representation to read
\be
\chi (a * b) = \chi (a) \chi (b), \qquad  \forall \; a, b \in {\cal R},
\ee
and set $c_1 = -1$, $c_2 = 0$ in formula (5.3).
(We can also adapt the dual point of view, def\/ining (right) representation
by the condition
\be
\chi (a * b) = - \chi (b) \chi (a),
\ee
and setting $c_1 = 0$, $c_2 = 1$ in formula (5.3). But we won't pursue
this avenue here, leaving it to the next Section.) Thus,
\be
\delta \psi (a,b) = \psi (a * b) - a . \psi (b),
\ee
where
\be
a . (\cdot) : = \chi (a) (\cdot).
\ee

All told, we def\/ine the coboundary operator $\delta: \ C^n \rightarrow
C^{n+1} $ by the formula
\be
\ba{l}
\ds \delta \psi (a_1,\ldots, a_n, a) = \sum^n_{i = 1} (-1)^{i+1}
[ \psi (\ldots \hat i \ldots, a_i * a) - a_i . \psi (\ldots\hat i
\ldots, a)]
\vspace{3mm}\\
\ds \qquad + \sum_{1 \leq i < j \leq n} (-1)^{i+j+1} \psi ([a_i, a_j]
\ldots \hat i \ldots \hat j \ldots, a).
\ea
\ee

For $n=0$, formula (5.10) is to be understood as
\be
\delta \psi (a_1) = - a_1 . \psi, \qquad  \psi \in {\cal M}.
\ee

The new extra sum in formula (5.10) doesn't destroy the property of
$\delta$ to preserve $\kappa$-skewsymmeetry.

\medskip

\noindent
{\bf Proposition 5.12.} $\delta^2 = 0.$

\medskip

\noindent
{\bf Proof.} We have seen above that $\delta^2 = 0$ on $C^0$, and
it's easy to verify that $\delta^2 = 0$ on $C^1$ and $C^2$. So let
$n \geq 3$.

\setcounter{equation}{12}
Setting
\be
\delta = \delta^{\mbox{\scriptsize old}} + \delta^{\mbox{\scriptsize new}},
\ee
where $\delta^{\mbox{\scriptsize old}}$ is given by formula (4.6),
and $\delta^{\mbox{\scriptsize new}}$ is given
by formula
\be
\delta^{\mbox{\scriptsize new}} \psi (a_1,\ldots, a_n, a) =
\sum^n_{i = 1} (-1)^i a_i.\psi (\ldots \hat i \ldots,a),
\ee
we have for $\nu = \delta \psi$:
\be
\ba{l}
\ds \nu (a_1, \ldots, a_n, z) = \sum^n_{i = 1} (-1)^{i+1} \psi (\hat i, iz)
\vspace{3mm}\\
\ds \qquad +\sum_{i < j} (-1)^{i+j+1} \psi ([i, j ] \hat i \hat j , z) +
\sum^n_{i = 1} (-1)^i a_i . \psi ( \hat i ,z),
\ea
\ee
\be
\ba{l}
\ds \delta \nu(y_1, \ldots, y_{n+1}, t) = \sum^{n+1}_{s=1} (-1)^{s+1} \nu
(\hat s, st)
 \vspace{3mm}\\
\ds \qquad  + \sum_{p<q} (-1)^{p+q+1} \nu([p, q] \hat p \hat q , t) +
 \sum^{n+1}_{\ell = 1} (-1)^\ell y_\ell . \nu (\hat \ell, t).
 \ea
 \ee

We shall work out separately each of the three sums in the expression
(5.16); since we have already verif\/ied in the preceding Section that
$(\delta^{\mbox{\scriptsize old}})^2 = 0$, we shall only keep track
of the extra terms coming out of the operator
$\delta^{\mbox{\scriptsize old}} \delta^{\mbox{\scriptsize new}} +
\delta^{\mbox{\scriptsize new}} \delta^{\mbox{\scriptsize old}} +
(\delta^{\mbox{\scriptsize new}} )^2$.

\noindent
({\it a}) \ We have:
\[
\nu ( \hat s, st) \doteq \sum_{i <s} (-1)^i y_i . \psi (\hat i\hat s, st)
 +\sum_{i > s} (-1)^{i+1} y_i . \psi ( \hat s\hat i , st).
\]
Multiplying this by $(-1)^{s+1}$ and summing on $s$, we get
\be
\{(5.16{\rm a})\} \doteq \sum_{a<b} (-1)^{a+b} \left[- y_a . \psi (\hat a
\hat b, b t) + y_b . \psi (\hat a \hat b , a t) \right];
\ee
({\it b}) \ We have:
\[
\ba{l}
\ds \nu ([p, q] \hat p \hat q, t) \doteq
- [y_p, y_q] . \psi (\hat p\hat q, t) + \sum_{\alpha< p} (-1)^{\alpha+1}
y_\alpha . \psi ([p, q] \hat \alpha \hat p \hat q  , t)
\vspace{3mm}\\
\ds \qquad + \sum_{p< \alpha<q} (-1)^\alpha y_\alpha . \psi ([p, q]
\hat p \hat \alpha \hat q, t) + \sum_{\alpha > q} y_\alpha . \psi ([p, q]
\hat p \hat q  \hat \alpha, t).
\ea
\]
Multiplying this by $(-1)^{p+q+1}$ and summing on $\{p<q\}$, we f\/ind:
\be
\{(5.16{\rm b})\} \doteq
\sum_{p<q} (-1)^{p+q} [y_p, y_q] . \psi ( \hat p \hat q, t)
\ee
\renewcommand{\theequation}{\arabic{section}.\arabic{equation}{\rm a}}
\setcounter{equation}{18}
\be
\qquad + \sum_{a<b<c} (-1)^{a+b+c} \{ y_a . \psi ([b,c] \hat a \hat b
\hat c,t)
\ee
\renewcommand{\theequation}{\arabic{section}.\arabic{equation}{\rm b}}
\setcounter{equation}{18}
\be
\qquad \qquad -y_b . \psi ([a, c] \hat a \hat b \hat c ,t)
\ee
\renewcommand{\theequation}{\arabic{section}.\arabic{equation}{\rm c}}
\setcounter{equation}{18}
\be
\qquad \qquad + y_c . \psi ([a, b] \hat a \hat b \hat c , t)\};
\ee
({\it c}) \ We have:
\[
\ba{l}
\ds y_\ell . \nu (\hat \ell  , t)= y_\ell .
\Biggl\{ \sum_{i < \ell} (-1)^{i+1} \psi (\hat i \hat \ell, it)
+ \sum _{i>\ell} (-1)^i \psi (\hat \ell \hat i, it)
\vspace{3mm}\\
\ds \qquad + \sum_{i<j<\ell} (-1)^{i+j+1} \psi([i, j]\hat i \hat j \hat \ell,
t) + \sum_{i < \ell < j} (-1)^{i+j} \psi ([i, j] \hat i \hat \ell\hat j , t)
\vspace{3mm}\\
\ds \qquad + \sum_{\ell < i < j} (-1)^{i+j+1} \psi ([i, j] \hat \ell \hat i
\hat j , t) +  \sum_{i<\ell} (-1)^i y_i . \psi ( \hat i \hat \ell , t)
 +\sum_{i > \ell} (-1)^{i+1} y_i . \psi ( \hat \ell \hat i ,
t)\Biggr\} .
\ea
\]
Multiplying all this by $(-1)^\ell$ and summing on $\ell$, we obtain:
\renewcommand{\theequation}{\arabic{section}.\arabic{equation}}
\setcounter{equation}{19}
\be
\{(5.16{\rm c})\}=
\sum_{a<b} (-1)^{a+b} \{ - y_b . \psi (\hat a\hat b, at) + y_a .
\psi ( \hat a \hat b, bt ) \}
\ee
\renewcommand{\theequation}{\arabic{section}.\arabic{equation}{\rm a}}
\setcounter{equation}{20}
\be
\qquad + \sum_{a<b<c} (-1)^{a+b+c} \{ -y_c . \psi ([a, b] \hat a \hat b
\hat c, t)
\ee
\renewcommand{\theequation}{\arabic{section}.\arabic{equation}{\rm b}}
\setcounter{equation}{20}
\be
\qquad \qquad +y_b . \psi ([a, c] \hat a \hat b \hat c , t)
\ee
\renewcommand{\theequation}{\arabic{section}.\arabic{equation}{\rm c}}
\setcounter{equation}{20}
\be
\qquad \qquad -y_a . \psi ([b, c] \hat a \hat b \hat c , t) \}
\ee
\renewcommand{\theequation}{\arabic{section}.\arabic{equation}}
\setcounter{equation}{21}
\be
\qquad + \sum_{a<b} (-1)^{a+b} \{ y_b . (y_a . \psi (\hat a\hat b,
t)) - y_a . (y_b . \psi (\hat a\hat b, t))\}.
\ee

The cancellation scheme is:

\begin{enumerate}
\item[1)] (5.17) and (5.20);

\item[2)] (5.18) and (5.22), since the action $\chi$ of ${\cal R}$ on
${\cal M}$ is a representation:
\be
\chi ([a, b]) = [ \chi (a), \chi (b)], \qquad \forall \; a, b \in
{\cal R};
\ee

\item[3)] (5.19) and (5.21).
\end{enumerate}

(Notice that $(\delta^{\mbox{\scriptsize new}})^2 \not= 0$.)
\hfill $\blacksquare$

\medskip

\noindent
{\bf Remark 5.24.} The coboundary operator $\delta$ (5.10) does
{\it not} reduce to the one of the Hochschild complex [2] when ${\cal R}$
is an associative algebra, even though the cochain spaces are identical
in both cases.

\medskip

\noindent
{\bf Remark 5.25.} When ${\cal M} = {\cal R}$ and the natural def\/inition
of representation is used, one arrives at a new complex by considering
deformations of the quasiassociative algebra ${\cal R}$, exactly like
the Hochschild complex on $C^\bullet({\cal R}, {\cal R})$
is arrived at in the associative case~[1]. This new complex is closely
related to the Hochschild one, and it is still dif\/ferent from the
one constructed above.

\setcounter{equation}{0}

\section{Dual Point of View, Homology}

The extended complex (5.10) of the preceding Section was based on the
notion of representation of a quasiassociative algebra ${\cal R}$
as a linear map $\chi: {\cal R} \rightarrow \mbox{End} ({\cal M})$
satisfying the condition
\be
\chi (a * b) = \chi (a) \chi (b).
\ee
There was a second version of representation, formula (5.7):
\be
\chi (a * b) = - \chi (b) \chi (a);
\ee
this choice was left unexamined. Let's examine it now.

These two choices lead to two dif\/ferent formulae for the coboundary
operator $\delta: C^1 \rightarrow C^2$,
\be
\delta \psi (a, b) = \psi (a * b) - \chi (a) \psi (b),
\ee
\be
\delta \psi (a, b) = \psi (a * b) + \chi (b) \psi (a).
\ee
The f\/irst direction was pursued in the preceding Section.
The second one, as is easy to discover by considering the hypothetical
map $\delta: C^2 \rightarrow C^3$, leads nowhere. Why is it so?

Let ${\cal M}^* = \mbox{Hom}_K ({\cal M}, K)$ be the dual space to
${\cal M}$. Since ${\cal R}$ acts on ${\cal M}$, it also acts on
${\cal M}^*$ in the dual way:
\be
\langle  \chi^d (a) (m^*), m\rangle  = - \langle m^*, \chi (a) (m)\rangle.
\ee
Hence,
\[
\ba{l}
\langle \chi^d (a * b) (m^*), m\rangle  = -
\langle  m^*, \chi (a * b) (m) \rangle
\vspace{2mm}\\
\ds \qquad = - \langle m^*,\chi (a) \chi (b) (m)
\rangle  =  - \langle \chi^d (b) \chi^d (a) (m^*), m\rangle,
\ea
\]
so that
\be
\chi^d (a * b) = - \chi^d (b) \chi^d (a).
\ee

Thus, our second version of representation, (6.2), is in fact dual to
the f\/irst one, (6.1). Therefore, this def\/inition is suited not for
cohomology but for the dual object, homology. Def\/ining the $n$-chains as
\be
C_0 = C_0({\cal R}, {\cal N}) = {\cal N}; \qquad
C_n = C_n ({\cal R}, {\cal N}) = {\cal N} \otimes {\cal R}^{\otimes
n}, \qquad n \in {\mathbf N},
\ee
where ${\cal N}$ is a ${\cal R}$-module on which ${\cal R}$ acts according
to formula (6.2):
\be
{\bar n}^\bullet (a * b) = - {(\bar n}^\bullet a)^\bullet b,
\qquad \bar n \in {\cal N}, \quad a, b \in {\cal R},
\ee
in the suggestive notation of the right action, we def\/ine the dif\/ferential
$\partial: C_n \rightarrow C_{n-1}$ by the rule:
\be
\ba{l}
\ds \partial (\bar n \otimes a_1 \otimes \ldots \otimes a_n) =
\sum^{n-1}_{i=1} (-1)^{i+1} \bar n \otimes \ldots \hat i
\ldots  a_i * a_n
\vspace{3mm}\\
\ds \qquad + \sum_{1 \leq i <j < n} (-1)^{i+j+1} \bar n \otimes
[a_i, a_j] \ldots  \hat i \ldots \hat j \ldots +
\sum^{n-1}_{i=1} (-1)^{i+1} (n^\bullet a_i) \otimes \ldots \hat i \ldots,
\vspace{3mm}\\
\ds \partial (C_0) = 0, \qquad  \partial (\bar n \otimes a) =
{\bar n}^\bullet a.
\ea
\ee

Since this formula satisf\/ies the duality relation
\be
\langle  \partial \Psi, \psi \rangle  = \langle \Psi, \delta \psi\rangle
\ee
for the case $\Psi \in C_n ({\cal R}, {\cal N}) \approx (C^n
({\cal R}, {\cal N}^*))^* $, $\psi\in C^{n-1} ({\cal R}, {\cal N}^*)$,
we have $\partial^2 = 0$ as a matter of course;
it is assumed that the chains considered are $\kappa$-skewsymmetric
for some $\kappa \geq 2$, exactly like the cochains.

\medskip

\noindent
{\bf Remark 6.12.} The Hochschild coboundary operator on $C^1 =\mbox{Hom}
({\cal R}, {\cal M}) $ acts by the rule:
\setcounter{equation}{12}
\be
\delta \psi (a_1, a_2) = a_1 . \psi (a_2) - \psi (a_1 a_2) +
\psi (a_1)^\bullet a_2,
\ee
where ${\cal R}$ is associative, ${\cal M}$ is an ${\cal R}$-bimodule,
and the right action of ${\cal R}$ on ${\cal M}$ is an anti-action from
the pont of view of our def\/inition~(6.2). We see that formulae~(6.3)
and~(6.4) each contribute about half to the Hochschild formula~(6.13).
There must be some underlying reason for such split.

\setcounter{equation}{0}

\section{Dif\/ferential Algebra Viewpoint}

Suppose our basic ring $K$ is a dif\/ferential ring, with a derivation
$\partial: K \rightarrow K$. Then the formula
\be
[X, Y] = XY' - X' Y, \qquad (\cdot)' = \partial (\cdot), \quad X, Y \in K,
\ee
makes $K$ into a Lie algebra ${\cal D}_1 = {\cal D}_1 (K)$,
the Lie algebra of vector f\/ields. The bilinear form $\omega$ on $K \times K$,
\be
\omega (X, Y) = XY'''
\ee
is skewsymmetric:
\be
\omega (X, Y) \sim - \omega (Y, X),
\ee
and is a generalized 2-cocycle on ${\cal D}_1$:
\be
\omega ([X, Y], Z) + \omega ([Y, Z], X) + \omega ([Z, X], Y) \sim 0,
\ee
where $(\cdot) \sim 0$ means that $(\cdot) \in \mbox{Im}\; \partial$.

When
\be
K = k \left[x, x^{-1}\right]
\ee
and
\be
\partial = d/dx,
\ee
$k$ being some number f\/ield or such, the Lie algebra ${\cal D}_1$ is
isomorphic to the centerless Virasoro algebra under identif\/ication
\be
e_n = x^{1-n} {d \over dx}, \qquad X = \sum_n X_n e_n.
\ee
As far as the Virasoro 2-cocycle is concerned, let
\be
\mbox{Res}: \;  k \left[x, x^{-1}\right] \rightarrow k
\ee
be the map isolating the $x^{-1}$-coef\/f\/icient, so that
\[
\mbox{Res} \circ \partial = 0.
\]
Then
\be
\ba{l}
\mbox{Res} (\omega (e_n, e_m)) = Res \left(x^{1-n} ( 1 - m) (-m) (-1-m) x^{-2-m}
\right)
\vspace{2mm}\\
\qquad = \delta^0_{n+m} (n+1) n (n-1) =
 \left(n^3 - n\right) \delta^0_{n+m}.
\ea
\ee

Below we construct a quasiassociative structure on
$K = k\left[x, x^{-1}\right]$ and the corresponding generalized
2-cocycle on it, so that formulae~(1.5)
and~(3.11) are recovered as localizations.

Let
\be
{\cal O} = x {d \over dx} -1
\ee
and set
\be
u * v = (1 - \epsilon {\cal O})^{-1} x^{-1} u (1 - \epsilon {\cal O})
{\cal O}(v),
\ee
\be
\hat \Omega (u, v) = x^{-3} {\cal O}^2 (1 + \epsilon {\cal O}) (u) \cdot v.
\ee
Since
\be
{\cal O} \left(x^{1-q}\right) = - q x^{1 - q},
\ee
we get
\be
\ba{l}
\ds x^{1 - p} * x^{1-q} = (1 - \epsilon {\cal O})^{-1}
\left(x^{-1} x^{1-p} (1 +\epsilon q) (-q) x^{1 - q}\right)
\vspace{3mm}\\
\ds \qquad = - q(1 - \epsilon q) (1 - \epsilon {\cal O})^{-1} x^{1-p-q} =
-{q(1 + \epsilon q) \over 1 + \epsilon (p+q)} x^{1 - p - q}.
\ea
\ee
This is formula (1.5). It implies that we have a correct quasiassociative
multiplication on $k\left[x, x^{-1}\right]$, with
\be
u*v - v*u = uv' - u' v.
\ee

The 2-cocycle story is more interesting. Recall how the notion of the
{\it generalized} 2-cocycle on a Lie algebra, equation~(7.4), appears:
from the classif\/ication of {\it affine} Hamiltonian operators,
with the linear part being attached to a Lie algebra, say ${\cal G}$,
and the constant part being a {\it generalized} 2-cocycle on this Lie
algebra~[4]. Aposteriori one can put all this into a variational complex
([4], p.~204) $\delta: \ \mbox{Dif\/f}\;(\wedge^n {\cal G}, K)_v \rightarrow
\ \mbox{Dif\/f}\; (\wedge^{n+1} {\cal G}, K)_v$, where subscript ``$v$''
signif\/ies that dif\/ferential forms dif\/fering by $Im\; \partial$
are to be identif\/ied; the generalized 2-cocycle condition~(7.4)
is then simply
\be
\delta \omega (X, Y, Z) \sim 0.
\ee

We shall now apply the same variational leap-forward to the complex
$C^\bullet ({\cal R}, K)$ of Section~5, considering cochains modulo
$\mbox{Im}\; \partial$.
A generalized 2-cocycle $\hat \Omega$ then satisf\/ies the dif\/ferential
version of the equality~(3.4):
\be
\delta \hat \Omega (u, v, w) = \hat \Omega (v, u*w) - \hat \Omega
(u, v*w) + \hat \Omega ([u, v], w) \sim 0.
\ee

It is unclear to me at the moment exactly what question such a variational
2-cocycle answers to, and repeated appeals to noncommutative dif\/ferential
geometry in the sense of Allan Connes haven't helped so far; nevertheless,
we have

\setcounter{equation}{19}
\medskip

\noindent
{\bf Proposition 7.19.} {\it (i)  Let $\hat \Omega$ be a generalized
2-cocycle on a differential quasiassociative algebra ${\cal R}$. Then
\be
\hat \omega (u, v) = \hat \Omega (u, v) - \hat \Omega (v, u)
\ee
is a generalized 2-cocycle on the Lie algebra $Lie ({\cal R})$;\\
(ii) The symplectic form on $T^*{\cal R}$ is a generalized 2-cocycle on
$T^*{\cal R}$.}

\medskip

\noindent
{\bf Proof.} (i) We have,
\be
\ba{l}
\ds \hat \omega ([u, v], w) = \hat \Omega ([u, v], w) - \hat \Omega (w, [u, v])
\vspace{3mm}\\
\ds \qquad   {\mathop {\sim}\limits^{\mbox{\scriptsize \rm [by (7.18)]}}} \ \
 \hat \Omega (u, v w) - \hat \Omega (v, uw) - \hat \Omega (w, uv - vu).
\ea
\ee
Hence,
\[
\ba{l}
\hat \omega ([u, v], w) + c.p. \sim (\hat \Omega (u, vw) + c.p.) -
(\hat \Omega (v, uw) + c.p.) - (\hat \Omega (w, uv - vu) + c.p.)
\vspace{3mm}\\
\ds \qquad
= (\hat \Omega (w, uv) + c.p.) - (\hat \Omega (w, vu) + c.p) - (\hat \Omega
(w, uv - vu) + c.p.) = 0;
\ea
\]

(ii)  By formula (1.1), $T^*{\cal R}$ has the multiplication
\renewcommand{\theequation}{\arabic{section}.\arabic{equation}{\rm a}}
\setcounter{equation}{21}
\be
\left(\begin{array}{c}u \\ \bar u \end{array} \right) *
\left(\begin{array}{c}v \\ \bar v \end{array} \right) =
\left(\begin{array}{c} u * v \\ u * \bar v \end{array} \right),
\qquad u, v \in {\cal R}, \quad \bar u, \bar v \in {\cal R}^*,
\ee
where
\renewcommand{\theequation}{\arabic{section}.\arabic{equation}{\rm b}}
\setcounter{equation}{21}
\be
\langle  u * \bar v, w \rangle  \sim - \langle  \bar v, u * w \rangle.
\ee
The symplectic form $\hat \Omega $ is
\[
\hat \Omega \left( \left(\begin{array}{c} u \\ \bar u \end{array} \right),
\left(\begin{array}{c} v \\ \bar v \end{array} \right) \right) = \langle
 \bar u, v\rangle - \langle  \bar v, u \rangle.
 \]
Hence,
\renewcommand{\theequation}{\arabic{section}.\arabic{equation}{\rm a}}
\setcounter{equation}{22}
\be
\ba{l}
\ds \hat \Omega \left( \left(\begin{array}{c} v \\ \bar v \end{array} \right),
\left(\begin{array}{c} u \\ \bar u \end{array} \right) *
\left(\begin{array}{c} w \\ \bar w \end{array}\right) \right) =
\langle \bar v, u * w \rangle  - \langle  u * \bar w, v\rangle
\vspace{3mm}\\
\ds \qquad \sim  \langle \bar v, u * w \rangle  + \langle \bar w, u * v
\rangle,
\ea
\ee
\renewcommand{\theequation}{\arabic{section}.\arabic{equation}{\rm b}}
\setcounter{equation}{22}
\be
- \hat \Omega \left( \left(\begin{array}{c} u \\ \bar u \end{array}\right),
\left(\begin{array}{c} v \\ \bar v \end{array} \right)*\left(\begin{array}{c}
w \\ \bar w \end{array} \right) \right)
\sim - \langle  \bar u, v * w \rangle  - \langle \bar w, v * u \rangle,
\ee
\renewcommand{\theequation}{\arabic{section}.\arabic{equation}{\rm c}}
\setcounter{equation}{22}
\be
\ba{l}
\ds \hat \Omega \left( \left[ \left(\begin{array}{c} u \\ \bar u \end{array}
\right), \left(\begin{array}{c} v \\ \bar v \end{array} \right) \right],
\left(\begin{array}{c} w\\  \bar w \end{array} \right) \right) =
\Omega \left( \left(\begin{array}{c} [u, v] \\ u * \bar v - v * \bar u
\end{array} \right), \left(\begin{array}{c} w \\ \bar w \end{array}
 \right) \right)
 \vspace{3mm}\\
 \ds \qquad = \langle  u * \bar v - v * \bar u , w \rangle  -
 \langle  \bar w, [ u, v] \rangle
 \vspace{3mm}\\
 \ds \qquad \sim \  -\langle \bar v, u * w \rangle  +
 \langle  \bar u, v * w \rangle  - \langle  \bar w, u * v - v * u\rangle.
 \ea
 \ee
Adding the expressions (7.23a-c) up, we get zero. \hfill $\blacksquare$

\medskip

Let us now verify that $\hat \Omega$ given by formula (7.13) is indeed a
generalized 2-cocycle. We have:
\renewcommand{\theequation}{\arabic{section}.\arabic{equation}}
\setcounter{equation}{23}
\be
\ba{ll}
1) &  \hat \Omega (v, u w) = x^{-3} {\cal O}^2 (1 +
\epsilon {\cal O}) (v) \cdot (1 - \epsilon {\cal O})^{-1} x^{-1} u
(1 - \epsilon {\cal O}) {\cal O}(w)
\vspace{2mm}\\
& \ds \qquad  \sim [1+ \epsilon ({\cal O} + 3)]^{-1} x^{-3} {\cal O}^2
(1 + \epsilon {\cal O}) (v) \cdot x^{-1} u (1 - \epsilon {\cal O}) {\cal O}
(w)
\vspace{3mm}\\
& \ds \qquad  \sim \left\{ - ({\cal O} + 3 ) [1 + \epsilon ({\cal O} + 3) ]
x^{-1} u [1 + \epsilon ({\cal O}+ 3)]^{-1} x^{-3} {\cal O}^2 (1 + \epsilon
{\cal O}) (v) \right\} \cdot w,
\ea \hspace{-16.8pt}
\ee
where we used the universal relation
\be
(1) A(2) \sim A^\dagger (1) \cdot (2)
\ee
for the adjoint operator, and the particular relation
\be
{\cal O}^\dagger = - ({\cal O} + 3)
\ee
for our operator $\ds {\cal O} = x {d \over dx} -1$.

Now, since
\be
({\cal O} + 3) x^{-3} = x^{-3} {\cal O},
\ee
formula (7.24) can be rewritten as
\renewcommand{\theequation}{\arabic{section}.\arabic{equation}{\rm a}}
\setcounter{equation}{27}
\be
\left\{ -x^{-4} ({\cal O}-1) [1+ \epsilon ({\cal O} -1)] u
{\cal O}^2 (v) \right\} \cdot w;
\ee
\renewcommand{\theequation}{\arabic{section}.\arabic{equation}{\rm b}}
\setcounter{equation}{27}
\be
\ba{ll}
2) & - \hat \Omega (u, vw) \sim \left\{ x^{-4} ({\cal O} - 1)
 [1 + \epsilon ({\cal O} -1)] v {\cal O}^2 (u) \right\} \cdot w;
\ea
\ee
\renewcommand{\theequation}{\arabic{section}.\arabic{equation}{\rm c}}
\setcounter{equation}{27}
\be
\ba{ll}
3) & \hat \Omega ([u, v], w ) = \left\{ x^{-3} {\cal O}^2 (1 + \epsilon
{\cal O}) (uv' -u' v) \right\} \cdot w.
\ea
\ee

Adding up the expressions (7.28a-c), we arrive at the equivalent relation
to be verif\/ied:
\renewcommand{\theequation}{\arabic{section}.\arabic{equation}}
\setcounter{equation}{28}
\be
x^{-3} {\cal O}^2 (1 + \epsilon {\cal O}) (uv' - u'v) = x^{-4}
({\cal O} -1)[1 + \epsilon ({\cal O} - 1)] \left[ v {\cal O}^2 (u) - u
{\cal O}^2 (v)\right].
\ee
Since
\be
x^{-1} ({\cal O} -1) = {\cal O} x^{-1},
\ee
equality (7.29) reduces to
\be
{\cal O} (uv' - u' v) = x^{-1} \left[u {\cal O}^2 (v) - v {\cal O}^2 (u)\right].
\ee
Now,
\be
{\cal O}^2 = x^2 {d^2 \over dx^2} - x {d \over dx} + 1,
\ee
so that
\[
\ba{l}
\ds x^{-1} \left[u {\cal O}^2 (v) - v {\cal O}^2 (w)\right] =
u (x v'' - v') -v (xu''- u')
\vspace{3mm}\\
\ds \qquad = x (uv'' - u'' v) - (u v' - u' v) =
\left(x {d \over dx} - 1\right) (uv' - u' v) = {\cal O} (uv' - u'v).
\ea
\]

It remains to perform the last step:  to calculate $\mbox{Res}\; \hat \Omega
\left(x^{1-p}, x^{1-q}\right)$ and to compare the result
with the formulae (3.11,12). We have:
\be
\hat \Omega \left(x^{1-p} , x^{1-q}\right) = x^{-3} {\cal O}^2
(1 + \epsilon {\cal O}) \left(x^{1-p} \right)
\cdot x^{1-q} \ \ {\mathop{=}\limits^{\mbox{\scriptsize \rm [by (7.14)]}}} \ \
(-p)^2 (1 - \epsilon p) x^{-1-p-q},
\ee
so that
\be
\mbox{Res}\; \hat \Omega \left(x^{1-p}, x^{1-q}\right) = p^2 (1 - \epsilon p)
\delta^0_{p+q} = - \epsilon \left(p^3 - \epsilon^{-1} p^2\right)
\delta^0_{p+q}.
\ee
We see that we have to multiply $\hat \Omega$ by
$- {1 \over 2} \epsilon^{-1}$, and also to add to it the trivial
2-cocycle proportional to the $*$ product. From formula~(7.15) we f\/ind:
\be
\mbox{Res} \; x^{-2} \left(x^{1-p} * x^{1-q} \right) = p (1 - \epsilon p)
\delta^0_{p+q} .
\ee
Thus, the correctly normalized generalized 2-cocycle has the form
\be
\hat \Omega^{\mbox{\scriptsize new}} (u, v) = - {1 \over 2}
\epsilon^{-1} x^{-3} {\cal O}^2
(1 + \epsilon {\cal O}) (u) \cdot v -
{1 \over 2} x^{-2} (1 - \epsilon {\cal O})^{-1} x^{-1} u
(1 - \epsilon {\cal O}){\cal O} (v).
\ee

\setcounter{equation}{36}

\noindent
{\bf Remark 7.36.} Consider the Lie algebra ${\cal D}_n = {\cal D}_n(K)$
``of vector f\/ields on ${\mathbf R}^n$'', with the commutator
\be
[X,Y]^i = \sum^n_{s=1} (X^s Y^i,_s - Y^s X^i,_s) , \qquad  X, Y \in K^n,
\ee
where
\be
(\cdot ),_s = \partial_s (\cdot) ,
\ee
and $\partial_1,\ldots, \partial_n: \ K \rightarrow K$ are $n$ commuting
derivations. Localizing $K$ as $k [x_1,\ldots, x_n, x_1^{-1}$, $\ldots,
x_n^{-1}]$ and  taking as the basis of ${\cal D}_n (K)$
\be
e^i_\sigma = x^{1 _{i}-\sigma} \partial_i = x_1^{-\sigma_{1}} \ldots
x_n^{- \sigma_n} x_i \partial_i, \qquad \sigma \in {\mathbf Z}^n, \quad
\partial_i = \partial/\partial x_{i },
\ee
we f\/ind the $n$-dimensional analog of the centerless Virasoro algebra:
\be
[e^i_\sigma, e^j_\nu] = (\delta_{ij} - \nu_i) e^j_{\sigma + \nu} -
(\delta_{ij} - \sigma_j) e^i_{\sigma + \nu}.
\ee
This Lie algebra does not seem to have a quasiassoactive representation of
the form~(1.5) for $n > 1$, but it does have a quasiassociative
representation generalizing formula~(2.11):
\be
e^i_\sigma * e^j_\nu = (\lambda \delta_{ij} - \nu_i) e^i_{\sigma + \nu},
\qquad \lambda = \mbox{const}.
\ee

\renewcommand{\theequation}{{\rm A}1.\arabic{equation}}
\setcounter{equation}{0}

\subsection*{Appendix 1. Virasoro Algebra Does Not Come from
an Associative One}

Suppose we have a ${\mathbf Z}$-graded multiplication on the basis
$\{e_p \mid p \in G$, a commutative ring$\}$, of the form
\be
e_i e_j = g (i, j) e_{i + j},
\ee
such that
\be
e_i e_j - e_j e_i = (i - j) e_{i+j}, \qquad \forall \; i, j \in G,
\ee
and
\be
(e_i e_j) e_k = e_i (e_j e_\kappa), \qquad \forall \; i, j, \kappa \in G.
\ee

Let us show that such representation is impossible.

We f\/irst rewrite the boundary condition (A1.2) as
\be
g (i, j) - g(j, i) = i-j.
\ee

Next, rewrite the associativity condition (A1.3) as
\be
g(i,j) g(i+j, \kappa) = g (j, \kappa) g (i, j+\kappa).
\ee

Now, set $j = \kappa = 0$ in formula (A1.5):
\be
g (i, 0) [g (i, 0) - g (0, 0)] = 0.
\ee

Further, set $j = i = 0$ in formula (A1.5):
\be
g (0, \kappa) [g (0, \kappa) - g (0, 0)] = 0.
\ee

Assume that $G$ has no zero divisors. From formula~(A1.6) we f\/ind:
\be
g (i, 0) = 0 \quad  \mbox{ or} \quad g(0,0),
\ee
while formula (A1.7) yields:
\be
g(0, \kappa) = 0 \quad \mbox{ or} \quad g(0, 0).
\ee
The last two equations contradict the boundary condition~(A1.4):
\[
g (r, 0) - g(0, r) = r.
\]

\renewcommand{\theequation}{{\rm A}2.\arabic{equation}}
\setcounter{equation}{0}

\subsection*{Appendix 2. Semidirect Sums of Quasiassociative Algebras}

Let ${\cal R}$ and ${\cal U}$ be quasiassociative algebras,
${\cal G} = Lie ({\cal R})$, ${\cal H} = Lie ({\cal U})$.
Let $\chi : {\cal G} \rightarrow \mbox{Der} ({\cal H})$ be a representation
of ${\cal G}$. The semidirect sum Lie algebra
${\cal G} {\mathop {\symb}\limits_{\chi}} {\cal H}$
is the vector space ${\cal G} \oplus {\cal H}$ with the commutator
\be
\left[ \left(\begin{array}{c} a \\ u \end{array} \right),
\left(\begin{array}{c} b \\ v \end{array}\right) \right] =
\left(\begin{array}{c} [a, b] \\ a . v - b . u + [ u, v] \end{array}
\right) , \qquad  a, b \in {\cal G}, \quad u, v \in {\cal H}.
\ee

Does the Lie algebra ${\cal G} {\mathop {\symb}\limits_{\chi}} {\cal H}$
have a quasiassociative representation?

\medskip

\setcounter{equation}{2}

\noindent
{\bf Proposition A2.2.} {\it Let $\chi: \  Lie({\cal R}) \rightarrow
\mbox{\rm Der} ({\cal U})$ be a representation. Define the semidirect sum
${\cal R} {\mathop {\symb}\limits_{\chi}} {\cal U}$
as the space ${\cal R} \oplus {\cal U}$ with the multiplication
\be
\left(\begin{array}{c} a \\ u \end{array} \right) *
\left(\begin{array}{c} b \\ b \end{array} \right)
= \left(\begin{array}{c} a * b \\ a . v + u * v \end{array} \right),
\qquad a, b \in {\cal R}, \quad u, v \in {\cal U}.
\ee
Then this multiplication is quasiassociative.}

\medskip

\noindent
{\bf Proof.} Dropping the $*$ notation for brevity, we have
\[
\hspace*{-14.2pt}
\left(\begin{array}{c}a \\ u \end{array} \right) \left( \left(\begin{array}{c}
b \\ v \end{array} \right) \left(\begin{array}{c} c \\ w \end{array}\right)
\right) = \left(\begin{array}{c} a \\ u \end{array}\right)
\left(\begin{array}{c} bc \\ b . w + v w \end{array} \right) =
\left(\begin{array}{c} a (bc) \\ a . (b . w + v w ) + u (b . w + vw)
\end{array} \right),
\]
\[
\left( \left(\begin{array}{c}a \\ u \end{array} \right)
\left(\begin{array}{c} b \\ v \end{array}\right)
\right) \left(\begin{array}{c}c \\ w \end{array} \right) =
\left(\begin{array}{c} ab \\ a . v + u v \end{array} \right)
\left(\begin{array}{c} c \\ w \end{array} \right) =
\left(\begin{array}{c} (a b) c \\ (a b) . w + (a . v + uv )w \end{array}
\right).
\]
Thus, we need to verify that
\[
\ba{l}
a . (b . w + vw) + u (b . w + vw) - (a b) . w - (a .
v + uv) w
\vspace{2mm}\\
 \qquad = b . (a . w + uw) + v (a . w + uw) - (b a) . w -
(b . u + vu) w.
\ea
\]
This can be rewritten as $0 {\mathop{=}\limits^{?}}$
\renewcommand{\theequation}{{\rm A}2.\arabic{equation}{\rm a}}
\setcounter{equation}{3}
\be
a . (b . w) - b . (a . w) - ((a b) . w - (b a). w)
\ee
\renewcommand{\theequation}{{\rm A}2.\arabic{equation}{\rm b}}
\setcounter{equation}{3}
\be
\qquad +a . (vw) - (a . v) w - v (a . w)
\ee
\renewcommand{\theequation}{{\rm A}2.\arabic{equation}{\rm c}}
\setcounter{equation}{3}
\be
\qquad + u (b. w) + (b. u) w - b. (u w)
\ee
\renewcommand{\theequation}{{\rm A}2.\arabic{equation}{\rm d}}
\setcounter{equation}{3}
\be
\qquad + u (vw) - (uv) w - v (u w) + (vu) w.
\ee
The f\/irst sum vanishes since $\chi$ is a representation of
$Lie \; ({\cal R})$;
the second and third sums vanish since
$Im (\chi) \subset \mbox{Der} ({\cal U})$; the fourth sum vanishes since
${\cal U}$ is quasiassociative. \hfill $\blacksquare$

\medskip

\noindent
{\bf Corollary A2.5.} {\it If ${\cal U}$ is abelian and $\chi: \
Lie ({\cal R}) \rightarrow \mbox{\rm End}({\cal U})$ is a representation,
then ${\cal R} {\mathop{\symb}\limits_{\chi}} {\cal U}$ is
quasiassociative.}

\medskip

\noindent
{\bf Proof.} $\mbox{Der} ({\cal U}) = \mbox{End} ({\cal U})$ for an
abelian ${\cal U}$. \hfill $\blacksquare$

\medskip

\noindent
{\bf Example A2.6.} Consider the Ehrenfest Lie algebra
${\cal G} (A)$, where $A$ is an arbitrary matrix, and the commutators
between basis elements are ([6], p.~274):
\renewcommand{\theequation}{{\rm A}2.\arabic{equation}}
\setcounter{equation}{6}
\be
[e_i, e_j] = [\bar e_i, \bar e_j] = 0, \qquad [e_i, \bar e_j] =A_{ji}\bar e_j.
\ee
In this case both ${\cal R}$ and ${\cal U}$ are vector spaces with trivial
multiplication,
\be
a * b = 0, \qquad  u * v = 0, \qquad  \forall \; a, b \in {\cal R},
\quad  u, v \in {\cal U},
\ee
and the representation $\chi$ acts by the fule
\be
e_i . \bar e_j = A_{ji} \bar e_j,
\ee
It {\it is} a representation of the abelian Lie algebra $Lie({\cal R})$,
since
\be
e_i . (e_j . \bar e_\kappa) = e_j . (e_i . \bar e_\kappa) =
A_{\kappa i}A_{\kappa j} \bar e_\kappa.
\ee
Hence, the Ehrenfest Lie algebra ${\cal G} (A)$ (A2.7) comes out of
the following quasiassociative multiplication:
\be
e_i e_j = \bar e_i \bar e_j = \bar e_i e_j = 0, \qquad
 e_i \bar e_j = A_{ji} \bar e_j.
\ee

\setcounter{equation}{12}

\noindent{\bf Remark A2.12.} Proposition A2.2 shows that
\be
Lie (R) {\mathop{\symb}\limits_{\chi}} Lie ({\cal U})
= Lie ({\cal R} {\mathop{\symb}\limits_{\chi}} {\cal U})
\ee
when ${\cal U}$ is abelian. Otherwise formula (A2.13) is
not necessarily true since $\mbox{Der} ({\cal U})$ is, in general,
smaller than $\mbox{Der} (Lie ({\cal U})):$

\medskip

\noindent
{\bf Proposition A2.14.} {\it (i)  $\mbox{\rm Der} ({\cal U})
\subset \mbox{\rm Der} (Lie ({\cal U}))$;\\
(ii)  If $\mbox{\rm Int}({\cal U}) \subset \mbox{\rm Der} ({\cal U})$
then ${\cal U}$ is associative. (Here $\mbox{\rm Int}({\cal U})$
denotes the space of maps $\{\mbox{\rm ad}_u: \ {\cal U} \rightarrow {\cal U} \ |u  \in
{\cal U} \}$.)}

\medskip

\noindent
{\bf Proof.} (i)  is well-known to be true for any algebra, not
necessarily associative or quasiassociative one;

(ii)  $\mbox{\rm ad}_u$ is a derivation of $Lie({\cal U})$ no matter whether ${\cal U}$
is quasiassociative or not. For $\mbox{ad}_u$ to be a derivation of ${\cal U}$,
we must have, for any $u, v, w \in {\cal U}$:
\[
0 = ad_u (vw) - (ad_u (v)) w -v ad_u (w)
= u (vw) - (vw) u - (uv - vu) w - v (uw - wu)
\]
\renewcommand{\theequation}{{\rm A}2.\arabic{equation}{\rm a}}
\setcounter{equation}{14}
\be
\qquad = u (v w) - (uv) w - v (uw) + (vu) w
\ee
\renewcommand{\theequation}{{\rm A}2.\arabic{equation}{\rm b}}
\setcounter{equation}{14}
\be
\qquad - (vw) u+ v (wu).
\ee
The f\/irst sum vanishes since ${\cal U}$ is quasiassociative.
The second sum vanishes if\/f ${\cal U}$ is associative.
\hfill $\blacksquare$

\renewcommand{\theequation}{{\rm A}3.\arabic{equation}}
\setcounter{equation}{0}

\subsection*{Appendix 3. Lie Algebras of Vector Fields on Lie Groups}

Formula
\be
X * Y = XY', \qquad X, Y \in C^\infty (S^{-1}), \quad {}' =
{d \over dz},
\ee
provides a quasiassociative structure on the Lie algebra of vector f\/ields
on the circle, ${\cal D}(S^1)$. Formula~[5]
\be
(X * Y)^i = \sum_s X^s Y^i,_{s}
\ee
provides a quasiassociative structure on the Lie algebra of vector f\/ields
on ${\mathbf R}^n, {\cal D} ({\mathbf R}^n)$. This suggests that
for some manifolds, similar structure exists for their Lie algebras of vector
f\/ields. (This will be proven below for $GL(n, {\mathbf R})$ and
$GL(n, {\mathbf C})$.) The parallelizable manifolds are the simplest,
and Lie groups are simpler still.

\medskip

\noindent
{\bf Proposition A3.3.} {\it Let ${\cal R}$ be a finite-dimensional
quasiassociative algebra over ${\mathbf R}$, ${\cal G} = Lie ({\cal R})$,
and $G$ a connected Lie group with the Lie algebra ${\cal G}$.
Then the Lie algebra of vector fields on $G$, ${\cal D} (G)$,
has a quasiassociative representation.}

\medskip
\setcounter{equation}{3}

\noindent
{\bf Proof.} Let $(e_i)$ be a basis in ${\cal R}$. Then
\be
e_i e_j = \sum_s \theta^s_{ij} e_s,
\ee
with some structure constants $\theta^s_{ij} \in {\mathbf R}$.
The quasiassociativity condition
\be
(e_i e_j) e_\kappa - e_i (e_j e_\kappa) = (e_j e_i) e_\kappa - e_j (e_i
e_\kappa), \qquad \forall \; i, j, \kappa,
\ee
translates into the equality
\be
\sum_s \left(\theta^s_{jk} \theta^r_{is} - \theta^s_{ij} \theta^r_{sk}\right) =
\sum_s \left(\theta^s_{ik} \theta^r_{js} - \theta^s_{ji} \theta^r_{sk}\right),
\ee
or
\be
\sum_s \left(\theta^s_{jk} \theta^r_{is} - \theta^s_{ik} \theta^r_{js}\right)=
\sum_s c^s_{ij} \theta^r_{sk},
\ee
where
\be
c^s_{ij} = \theta^s_{ij} - \theta^s_{ji}
\ee
are the structure constants of the Lie algebra ${\cal G} = Lie ({\cal R})$:
\be
[e_i, e_j] = e_i e_j - e_j e_i = \sum_s c^s_{ij} e_s.
\ee

Denote by $\hat e_i$ the left-invariant vector f\/ields on $G$ generated
by the elements $e_i \in {\cal G},$ so that
\be
\hat e_i \hat e_j - \hat e_j \hat e_i = \sum_j c^s_{ij} \hat e_s.
\ee
In this basis, every vector f\/ield on $G$ can be identif\/ied with a vector
from $C^\infty (G)^{\mbox{\scriptsize \rm dim}(G)}$:
\be
X \in {\cal D} (G) \  \Rightarrow
 \ X = \sum_i X^i \hat e_i, \qquad  X^i \in C^\infty (G).
\ee

For $X = \sum X^i \hat e_i$, $Y = \sum Y^j \hat e_j  \in {\cal D} (G)$,
set
\be
(X * Y)^r = \sum_\alpha X^\alpha \hat e_\alpha (Y^r) +
\sum_{\alpha \beta} X^\alpha Y^\beta \theta^r_{\alpha \beta}.
\ee

We are going to show that this multiplication makes ${\cal D} (G)
\approx C^\infty (G)^{\mbox{\scriptsize \rm dim}(G)}$ into a
quasiassociative algebra; the boundary conditions are
satisf\/ied since
\be
\ba{l}
\ds \sum_r (X * Y - Y * X)^r \hat e_r = \sum_{\alpha r} \left[X^\alpha
\hat e_\alpha (Y^r) - Y^\alpha \hat e_\alpha (X^r) \right] \hat e_r
\vspace{3mm}\\
\ds \qquad + \sum_{\alpha \beta r} X^\alpha Y^\beta c^r_{\alpha \beta}
\hat e_r =  \sum_{\alpha \beta} \left[X^\alpha \hat e_\alpha, Y^\beta \hat
e_\beta\right] =[X, Y].
\ea
\ee

Now,
\be
\ba{l}
\ds (X (YZ))^r =
\sum_\alpha X^\alpha \hat e_\alpha ((YZ)^r) + \sum_{\alpha
\beta} \theta^r_{\alpha \beta} X^\alpha (YZ)^\beta
\vspace{3mm}\\
\ds \qquad
= \sum_\alpha X^\alpha \hat e_\alpha \left(\sum_\mu Y^\mu \hat e_\mu
(Z^r)+ \sum_{\mu \nu} \theta^r_{\mu \nu} Y^\mu Z^\nu\right)
\vspace{3mm}\\
\ds \qquad + \sum_{\alpha \beta} \theta^r_{\alpha \beta} X^\alpha
\left(\sum_\mu Y^\mu \hat e_\mu (Z^\beta)
+ \sum_{\mu \nu} \theta^\beta_{\mu \nu} Y^\mu Z^\nu\right)
\vspace{3mm}\\
\ds \qquad =
\sum_{\alpha \mu} X^\alpha \hat e_\alpha (Y^\mu) \hat e_\mu (Z^r)
+ \sum_{\alpha \mu} X^\alpha Y^\mu \hat e_\alpha \hat e_\mu (Z^r)
\vspace{3mm}\\
\ds \qquad + \sum_{\alpha \mu \nu} X^\alpha \theta^r_{\mu \nu}
\left(\hat e_\alpha (Y^\mu) Z^\nu + Y^\mu \hat e_\alpha (Z^\nu)\right)
\vspace{3mm}\\
\ds \qquad +\sum_{\alpha \mu \nu} \theta_{\alpha \nu}^r X^\alpha Y^\mu \hat e_\mu
(Z^\nu) + \sum_{\alpha s \mu \nu} \theta^r_{\alpha s} \theta^s_{\mu \nu}
X^\alpha Y^\mu Z^\nu,
\ea
\ee
\be
\ba{l}
\ds ((X Y)Z)^r = \sum_\mu (XY)^\mu \hat e_\mu (Z^r) +
\sum_{s \nu} \theta^r_{s \nu} (XY)^s Z^\nu =
 \sum_{\mu \alpha} X^\alpha \hat e_\alpha (Y^\mu) \hat e_\mu (Z^r)
 \vspace{3mm}\\
\ds \qquad +\sum_{\mu \alpha \beta} \theta^\mu_{\alpha \beta} X^\alpha
Y^\beta \hat e_\mu (Z^r) + \sum_{s \nu} \theta^r_{s \nu}
\left(\sum_\alpha X^\alpha \hat e_\alpha (Y^s)+
\sum_{\alpha \beta} \theta^s_{\alpha \beta} X^\alpha Y^\beta \right) Z^\nu.
\ea
\ee
Thus,
\[
(X(YZ) - (XY) Z)^r
\]
\renewcommand{\theequation}{{\rm A}3.\arabic{equation}{\rm a}}
\setcounter{equation}{15}
\be
\qquad = \sum_{\alpha \mu} X^\alpha Y^\mu \hat e_\alpha \hat e_\mu
(Z^r)
\ee
\renewcommand{\theequation}{{\rm A}3.\arabic{equation}{\rm b}}
\setcounter{equation}{15}
\be
\qquad + \sum_{\alpha \mu \nu} X^\alpha Y^\mu \theta^r_{\mu \nu}
\hat e_\alpha (Z^\nu) + \sum_{\alpha \mu \nu} X^\alpha Y^\mu
\theta^r_{\alpha \nu} \hat e_\mu (Z^\nu)
\ee
\renewcommand{\theequation}{{\rm A}3.\arabic{equation}{\rm c}}
\setcounter{equation}{15}
\be
\qquad + \sum_{\alpha s \mu \nu} X^\alpha Y^\mu Z^\nu \left(\theta^r_{\alpha s}
\theta^s_{\mu \nu} - \theta^r_{s \nu} \theta^s_{\alpha \mu}\right)
\ee
\renewcommand{\theequation}{{\rm A}3.\arabic{equation}{\rm d}}
\setcounter{equation}{15}
\be
\qquad - \sum_{\mu \alpha \beta} \theta^\mu_{\alpha \beta} X^\alpha Y^\beta
\hat e_\mu (Z^r).
\ee
Interchanging $X$ and $Y$, subtracting the resulting expressions,
noticing that (A3.16b) is symmpletric in $(X, Y)$, and
using formulae (A3.10,8), we arrive at the following identity to be
verif\/ied: $0 {\mathop{=}\limits^{?}}$
\renewcommand{\theequation}{{\rm A}3.\arabic{equation}{\rm a}}
\setcounter{equation}{16}
\be
=\sum_{\alpha \mu s} X^\alpha Y^\mu \left(\theta^s_{\alpha \mu} -
\theta^s_{\mu\alpha}\right) \hat e_s (Z^r)
\ee
\renewcommand{\theequation}{{\rm A}3.\arabic{equation}{\rm b}}
\setcounter{equation}{16}
\be
+ \sum_{\alpha s \mu \nu} X^\alpha Y^\mu Z^\nu \left(\theta^r_{\alpha s}
\theta^s_{\mu\nu} - \theta^v_{s \nu} \theta^s_{\alpha \mu} -
\theta^r_{\mu s} \theta^s_{\alpha \nu} + \theta^r_{s \nu} \theta^s_{\mu \alpha}
\right)
\ee
\renewcommand{\theequation}{{\rm A}3.\arabic{equation}{\rm c}}
\setcounter{equation}{16}
\be
- \sum_{\mu \alpha \beta} \left(\theta^\mu_{\alpha \beta} -
\theta^\mu_{\beta\alpha} \right) X^\alpha Y^\beta \hat e_\mu (Z^r).
\ee

The expressions (A3.17a) and (A3.17c) cancel each other out.
The sum (A3.17b) vanishes due to the quasiassociativity condition (A3.6).
\hfill $\blacksquare$

\label{kupershmidt_4-lp}

\end{document}